\documentclass{amsart}
\usepackage{amssymb}
\usepackage{amsfonts}
\usepackage{amsbsy}

\usepackage{latexsym}
\usepackage[mathscr]{eucal}
\usepackage{verbatim}

\title[Commutators and linear spans of projections]{Commutators and
linear spans of projections in certain finite  C*-algebras}

\author{Victor Kaftal}

\address{Department of Mathematics\\
University of Cincinnati\\
P. O. Box 210025\\
Cincinnati, OH\\
45221-0025\\
USA}

\email{victor.kaftal@uc.edu}

\author{P. W. Ng}

\address{Department of Mathematics\\
University of Louisiana\\
217 Maxim D. Doucet Hall\\
P.O. Box 41010\\
Lafayette, Louisiana\\
70504-1010\\
USA}

\email{png@louisiana.edu}

\author{Shuang Zhang}

\address{Department of Mathematics\\
University of Cincinnati\\
P.O. Box 210025\\
Cincinnati, OH\\
45221-0025\\
USA}

\email{shuang.zhang@uc.edu}

\date{\today\\AMS Mathematics subject classification 46L35}

\theoremstyle{plain}
\newtheorem{theorem}{Theorem}[section]
\newtheorem{corollary}[theorem]{Corollary}
\newtheorem{lemma}[theorem]{Lemma}
\newtheorem{proposition}[theorem]{Proposition}
\newtheorem{remark}[theorem]{Remark}
\newtheorem{definition}[theorem]{Definition}

\newcommand{\be}{\begin{equation}\label}
\newcommand{\ee}{\end{equation}}
\newcommand{\bq}{\begin{equation*}}
\newcommand{\eq}{\end{equation*}}
\newcommand{\ba}{\begin{align*}}
\newcommand{\ea}{\end{align*}}
\newcommand{\bp}{\begin{proof}}
\newcommand{\ep}{\end{proof}}
\newcommand{\bL}{\begin{lemma}\label}
\newcommand{\eL}{\end{lemma}}
\newcommand{\bP}{\begin{proposition}\label}
\newcommand{\eP}{\end{proposition}}
\newcommand{\bC}{\begin{corollary}\label}
\newcommand{\eC}{\end{corollary}}
\newcommand{\bT}{\begin{theorem}\label}
\newcommand{\eT}{\end{theorem}}
\newcommand{\bR}{\begin{remark}\label}
\newcommand{\eR}{\end{remark}}
\newcommand{\bD}{\begin{definition}\label}
\newcommand{\eD}{\end{definition}}

\newcommand{\A}{\mathcal{A}}

\newcommand{\C}{\mathcal{C}}

\newcommand{\Mul}{\mathcal{M}}
\newcommand{\M}{\mathbb{M}}

\newcommand{\K}{\mathcal{K}}

\DeclareMathOperator{\TA} {T(\A)}
 \DeclareMathOperator{\QA} {QT(\A)}
\DeclareMathOperator{\QAff} {Aff\QA}
\DeclareMathOperator{\Aff}{Aff\TA}
 \DeclareMathOperator{\Ext}{Ext(\TA)}
  \DeclareMathOperator{\ext}{Ext(}
 
  \DeclareMathOperator{\her} {her}
 \DeclareMathOperator{\tr} {Tr}
\DeclareMathOperator{\T} {T}
\numberwithin{equation}{section}

\def\sideremark#1{\ifvmode\leavevmode\fi\vadjust{\vbox to0pt{\vss
\hbox to 0pt{\hskip\hsize\hskip1em
\vbox{\hsize2cm\tiny\raggedright\pretolerance10000
\noindent#1\hfill}\hss}\vbox to8pt{\vfil}\vss}}}

\begin{document}

\begin{abstract}
  Assume that  $\A$ is a  unital  separable simple C*-algebra
with real rank zero, stable rank one, strict comparison of
projections, and that its tracial simplex $\TA$ has a finite number
of extremal points.   We prove that every self-adjoint element $a$
in $\A$ with $\tau(a)=0$ for all $\tau\in \TA$ is the sum of two
commutators in $\A$ and that  that every positive element of $\A$ is
a linear combination of projections with positive coefficients.
Assume that $\A$ is as above
 but $\sigma-$unital. Then
an element (resp. a positive element) $a$ of $\A$   is a linear combination (resp. a
linear combination with positive coefficients) of
projections  if and only if $\bar\tau(R_a)<\infty$ for every $\tau\in
\TA$, and if and only if  , where $\bar\tau$ denotes the  extension of
$\tau$ to a tracial weight on $\A^{**}$  and $R_a\in \A^{**}$ denotes the range
projection of $a$. Assume that $\A$ is unital and as above but $\TA$ has infinitely many extremal points. Then  $\A$ is not the linear span of its projections. This result settles  two open problems of Marcoux in \cite {MarcouxIrishSurvey}.
\end{abstract}

\maketitle

\section{Introduction}\label{S:intro}
In the history of Operator Theory and Operator Algebras the study of
how bounded operators are composed of the fundamental building
blocks, projections, has attracted many researchers' attention.

In 1967 Fillmore \cite{Fillmore} found that every bounded
operator on a separable Hilbert space is made up of a linear
combination of 257 projections. Soon after, the number of needed
projections was reduced to 16 by Pearcy and Topping \cite
{PearcyToppingIdempotents} via Brown and Pearcy's characterization of commutators
\cite {BrownPearcyTypeI} and more recently to 10 by
Matsumoto \cite {Matsumoto}.

Pearcy and Topping (\cite {PearcyToppingIdempotents} and \cite {PearcyTopping}) proved that every element in properly infinite von
Neumann algebras or in certain type II$_1$ factors (Wright factors)
can be decomposed as  a (finite) linear combination of projections.
The same result was proven for the harder case of all
  type II$_1$ von Neumann algebras  by Fack and De La Harpe
\cite {FackDelaHarpe} and  then Goldstein and Paszkiewicz \cite {GoldsteinPaszkiewicz}
proved that the same holds  if and only if the von Neumann algebra does
not have a finite type I direct summand with infinite dimensional
center.

Of course the same conclusion cannot hold for all C*-algebras given
the lack of
  projections in some algebras (e.g., see Blackadar \cite {BlackadarProjectionless}). Thus, the following
   question arises naturally:\\
 (ALP):  Which C*-algebras are the (algebraic) linear span of projections?
 And if an algebra is not the  linear span of its projections, how to  characterize elements that are linear combination of projections?

It is natural to focus first on the algebras with
 the highest level of abundance of projections, namely, C*-algebras of real rank zero
    (\cite {BrownPedersen}). Every self-adjoint element in a C*-algebra of real rank zero  can be approximated by linear
  combinations of mutually orthogonal projections, namely  elements with finite spectrum.

Marcoux  has proved (\cite {MarcouxIndianaSpanProjections} and
\cite{MarcouxSmallNumberCommutators}, see also his survey \cite {MarcouxIrishSurvey}) that the following simple
 unital  C*-algebras   are the linear span of  projections. 

\begin{itemize}
\item Simple purely infinite ones.
\item AF-algebras with finitely many extremal tracial states.
\item AT-algebra with real rank zero and finitely many extremal tracial states.
\item Certain AH-algebras with real rank zero, bounded dimension growth, and finitely many extremal tracial states.
\end{itemize}

Related to the  (ALP) question is the following non-trivial question:\\

(PCP):  For which C*-algebras are all \textit {positive} elements
linear combinations of projections with positive coefficients?
(\textit {positive combinations of projections} for short). And if
not all, how to characterize   positive elements that are positive
combination of projections?

\

Fillmore observed (\cite {Fillmore})  that  a positive infinite
rank compact operator in $B(H)$ cannot be a positive combination of
projections; Fong and Murphy proved (\cite {FongMurphy}) that these
are the only bounded operators which are not.

 An analogous result holds in  $\sigma-$finite type II$_\infty$ von Neumann   factors where we proved in \cite [Corollary 3.5]{KNZFiniteSumsVNA}
  that all
   positive elements  are positive combinations of projections
   except those that have infinite range projection and belong to  the Breuer ideal generated by all finite projections.
    Moreover, all  positive elements in a  von Neumann factor of type  I$_n$, II$_1$, or $\sigma-$finite type III
   are  positive
     combinations     of projections (\cite [Theorem 2.12]{KNZFiniteSumsVNA}). In the non  $\sigma-$finite case or in
     von Neumann algebras with a nontrivial center,  a necessary and sufficient condition for a positive
     element to be a positive combination of projections is given in terms of central ideals and the central
      essential spectrum  (\cite [Theorem 2.12]{KNZFiniteSumsVNA}).

As usual, the purely infinite case is more tractable. We did prove
in \cite{KNZPISpan} that all positive elements in purely infinite
simple C*-algebras or in their multiplier algebras  are positive
linear combinations of projections. Note that all purely infinite
simple C*-algebras have  real rank zero
(\cite{ZhangWeylVonNeumann}). Finite C*-algebras of real rank zero,
however, are considerably harder.

In this article our targets are simple separable C*-algebras with real rank zero,
 stable rank one, which have  strict comparison of projections
 and have finitely many extremal tracial states. It may be somewhat surprising that  the determining factor turns out
to be  the number of  extreme points in $\TA$, which is a w*-compact simplex. We will show  that if the number of extreme
points is infinite then there are   positive elements  that are
not linear combinations of projections (Proposition \ref {P: no lin}); this settles two questions by Marcoux \cite
{MarcouxIrishSurvey}. If the number of  extreme points is finite, then every element  in the algebra is a linear combination of projections (Theorem \ref {T:lin comb}) and every positive element is a positive linear combination of projections (Corollary \ref {C:unital}).  This subsumes the above mentioned results for finite algebras by Marcoux. For non-unital algebras we will provide in Theorem \ref {T:main} and Corollary \ref{C:non pos} necessary and sufficient conditions for an element (resp. a positive element) to be a linear combination (resp. a positive combination ) of projections. This shows that neither the ALP property nor the PCP property are invariant under Morita equivalence.

 Key ingredients in our proofs are
 \begin{itemize}\item Embedding in
 $ \mathscr A$ a  unital simple AH-algebra $\mathscr C$ with real rank zero and  dimension
 growth bounded by 3 and having the same K-theory invariants,  based on a result of Lin (\cite
 {LinRR0AHEmbedding})
  and the  Elliott-Gong classification (\cite {ElliottGongRR0});
 \item Extending the construction
 of  Fack (\cite {Fack}) and Thomsen (\cite{ThomsenCommutator})  via certain inductive limit
  to  approximate  elements with zero traces  by a bounded number of commutators;
 \item  Marcoux's technique (\cite {MarcouxSmallNumberCommutators}) to 
 express  those element  as sum of commutators first,  and then as linear combination of projections;
 \item Brown's interpolation property (\cite{BrownIMP});
  \item the extension of traces on $\A$ to tracial weights on $\A^{**}$  of Combes \cite {Combes} and of Ortega, Rordam, and Thiel \cite{OrtegaRordamThiel}
    .
 \end{itemize}

\section{Preliminaries}\label{S:prel}
\subsection{The tracial simplex} \label{s:Tracial states}

\smallskip

If $\A$ is a unital simple C*-algebra of real rank zero, denote by $\TA$ the collection of the tracial states on $\A$.  It is well known that $\TA$ is a w*-closed convex subset of the state space  of  $\A$ and hence is w*-compact. $\TA$  is in fact a Choquet simplex \cite [Theorem 3.1.18] {Sakai}.

\smallskip

When $\A$ is a $\sigma$-unital but not unital simple C*-algebra of real rank zero, denote by $\widetilde \T (\A)$
 the collection of all nonzero lower semi-continuous, semifinite tracial weights on $\A $. Recall that a weight on $\A$ is
 a map $\tau : \A  _+ \longrightarrow [0,\infty ]$ such that
\begin{itemize}
\item $\tau (\lambda x ) =\lambda \tau (x) $ for all $x\in \A _+$ and
all $\lambda \ge 0$.

\item $\tau (x+y) = \tau (x) + \tau (y)$ for all $x, y\in \A  _+$.

A weight is called tracial if furthermore

\item $\tau (xx^*)=\tau (x^*x)$  for all $x\in \A $;

it is called semifinite (or densely defined) if

\item $\{a\in \A \mid \tau(a^*a)< \infty\}$ is dense in $\A$;

and it is called faithful if

\item $\tau(a)>0$  for all $0\ne a\in A _+$.
\end{itemize}
In the literature, tracial weights are also called  traces or extended traces, while tracial states are also called normalized traces.

Notice that lower-semicontinuous weights on a   C*-algebra of real
rank zero are completely determined by their values on the
projections of the algebra. Indeed,
  every positive element $a\in\A _+$  is the norm limit of a sequence of positive operators $b_n$ with finite spectrum
  (\cite{BrownPedersen}); in addition,
  the operators $b_n$ can be chosen such that $b_n\le a$, (e.g., see\cite [Lemma 2.3] {KNZ Mult}). Thus two
weights that agree on projections must agree on positive operators
with finite spectrum, and hence  by their lower-semicontinuity,  must agree also on all positive operators.

Furthermore, if $p\in \A$ is a nonzero projection,  by the
simplicity of $\A$ it follows that for every projection $q\in \A$
there is an $n\in \mathbb N$ for which $[q]\le n[p]$, where $[p]$
denotes the Murray-von Neumann equivalence class  of the projection
$p$ (e.g., see \cite [Corollary 6.3.6] {BlackadarBook}.) Thus if a
lower semincontinuous tracial weight $\tau $ is finite (resp.
nonzero) for one nonzero projection, then it must be finite (resp.
nonzero) for all nonzero projections and hence for all positive
operators with finite spectrum. Thus if $\tau$ is nonzero, then it
is faithful. Moreover, $\tau$ is semifinite if and only if $\tau(p)<
\infty$ for some non-zero projection.

Notice  further that every $\tau\in \widetilde \T(\A\otimes \K)$ is
uniquely determined by its values on $\A _+\otimes e_{11}$ where
$\{e_{ij}\}$ denotes the system of matrix units of the C*-algebra of
all compact operators $\K$ on a separable Hilbert space. It is well
known that every projection $p\in \A\otimes \K$ is unitarily
equivalent to a projection in $\M_n(\A)= \A \otimes \M_n(\mathbb C)$
for some $n\in \mathbb N$ and that any tracial weight on $(\A
\otimes \M_n(\mathbb C)) _+$ is determined by its values on $\A
_+\otimes e_{11}$. Thus if $\tau\in \widetilde \T(\A)$ and  $ \tr$
denotes the standard trace on $\K$, then $\tau \otimes \tr $ is the
unique extension of $\tau$ to an element of $\widetilde \T(\A\otimes
\K)$. This permits us to identify $ \widetilde \T(\A)$ with
$\widetilde \T(\A  \otimes \K).$

It follows from the  Brown's Stabilization Theorem \cite [Corollary
2.6] {BrownStable} that
  for every full projection $p\in \Mul(\A )$  there is an isometry $w\in \Mul(\A\otimes \K)$
  such that $$w(\A\otimes \K)w^*= p \A p \otimes \K .$$  Thus if $\tau\in  \widetilde \T(p\A p \otimes \K)$,
  then $\tau'(\cdot ):= \tau(w^*\cdot w)$ is the unique extension of $\tau$ to an element of $  \widetilde T(\A  \otimes \K)$.

Combining the  two observations above, every tracial state $\tau$ in
$\T(p\A p)$ uniquely extends to $\tau \otimes \tr \in
\widetilde\T(p\A p \otimes \K)$, and in turn to to $(\tau \otimes \tr)
(w^* \cdot w) \in  \widetilde \T(\A \otimes \K)$, and the restriction
of $(\tau \otimes \tr) (w^* \cdot w)$ to $\A _+\otimes e_{11}\cong \A
_+$ is therefore the unique extension of $\tau$ to an element of
$\widetilde \T(\A)$. Therefore, we can identify $\T(p\A p)$ with the
scaled space of $\widetilde \T(A)$, that is
\begin{equation} \label{e: trace space} \T(p\A p) = \big \{\frac {\tau}{\tau(p)} \mid \tau\in \widetilde \T(\A)\big\}.
\end{equation}

The next lemma  shows that the tracial simplex $\T(p\A p)$ does not depend on the choice of the nonzero projection $p$.

\bL{L:trace space} Let $p, q\in \A$ be nonzero projections. Then
there is a bijection $\psi$ between the sets of extremal points
$\ext T(p\A p))$ and $\ext T(q\A q))$ of $T(p\A p)$ and of $ T(q\A
q)$, respectively. Furthermore, if  $\ext T(p\A p)))$ is a finite
set, then the bijection $\psi$ extends to an affine homeomorphism
between $T(p\A p)$ and $T(q\A q)$. \eL

\begin{proof}

Define $\psi : \ext T(p\A p)) \longrightarrow \ext T(q\A q))$ by
 $ \psi \big(\frac {\tau } {\tau (p)}\big) =  \frac {\tau } {\tau (q)}.$
Then $\psi$ is a bijection. In fact, given $\frac { \tau }{\tau (p)}$
 an extreme point of $T(p\A p)$ for some $\tau\in \widetilde \T(\A)$,
 we show  that $\frac{\tau}{\tau (q)}$ is an extreme point of $T(q\A q)$.

We reason this by contradiction. Assume that the image $\psi
\big(\frac {\tau } {\tau (p)}\big) =\frac {\tau } {\tau (q)}$ is
written as the convex combination in $\T(q\A q)$
$$ \frac {\tau } {\tau (q)}= t \frac {\tau_1 } {\tau_1 (q)} +(1-t)
\frac {\tau_2 } {\tau_2(q)} \ \ \text{for some }\ \ 0<t<1.$$ Then
$$\frac {t \,\tau (q)\tau _1(p)}{\tau (p)
 \tau_1 (q)} + \frac {(1-t)\tau (q)  \tau_2 (p)} {\tau (p)\tau_2(q)}
 =1$$
and hence
\begin{align*}
 \frac {\tau } {\tau (p)}  &=  \frac { t\,\tau (q)} {\tau (p)\tau_1 (q)} \,\tau _1 + \frac {(1-t)\tau (q)  } {\tau (p)\tau_2(q)} \, \tau_2
\\&= \frac {t\, \tau (q)\tau _1(p)}{\tau (p)
 \tau_1 (q)} \frac {\tau _1}{\tau_1 (p)}+  \frac {(1-t)\tau (q)  \tau_2 (p)} {\tau (p)\tau_2(q)}\frac {\tau_2}{\tau _2 (p)}
 \end{align*}
would be a proper convex combination of tracial states   $\frac
{\tau _1}{\tau _1 (p)}$ and $\frac {\tau _2}{\tau _2 (p)}$ of $T(p\A
p)$, a contradiction. Thus $\psi $ maps extreme points to extreme
points. By switching the roles of $p$ and $q$, we see that  $\psi $
is one-to-one and onto.

If furthermore  $T(p\A p)$ has only finitely many extremal points $\{\tau_j\}_1^m$,
we can uniquely extend $\psi$  to an affine  map of the simplex $T(p\A p)$ onto the simplex
$T(q\A q)$ by setting $$\psi(\sum_{i=1}^m \lambda _i \tau_i) =\sum_{i=1}^m \lambda
_i \psi (\tau_i) \ \ \textit{for all} \ \
 0\le \lambda _i\le1 \ \ \text{with}\ \ \sum_{i=1}^m \lambda _i =1.$$
The  continuity of $\psi$ is then obvious.
\end{proof}
From now on    we will frequently identify $T(p\A p)$ with $\tilde
\T(\A)$ for any nonzero projection $p\in \A$, and denote both by
$\TA$ and use $\tau$ to denote both a
 tracial state on $p\A p$ as well as its extension to a lower semincontinuous semifinite tracial
 weight  on $\A _+$ or on $(\A \otimes \K)_+$.

\subsection{Continuous affine function on the tracial simplex}\label{s:aff}

 Recall that $\TA$ is a compact convex space. Let $\Aff$ be the space of all real-valued,
continuous, affine functions on $\TA$, equipped with the uniform
norm. $ \Aff$ is a closed subspace of the space of all real-valued
continuous functions on $\TA$, and hence is a (real) Banach space.

For every projection  $p\in \A\otimes \K$  let $\hat p$ denote the  evaluation map
$$\TA\ni \tau\to \hat p(\tau):=(\tau\otimes \tr )(p).  $$
It is elementary to see that $\hat p\in \Aff$.

Notice that if $p\in \A$, then of course $ \hat p(\tau)=\tau(p)$. Notice also that if
 $\A$ is unital and $a=a^*\in \A$, then the evaluation map $$\TA\ni \tau\to \hat a(\tau):=\tau(a)
  $$ also belongs to  $\Aff$.

Now  consider the map
\be{e:isom}
\Phi:  \Aff \ni f \to \Phi(f)=\{f(\tau)\}_{\tau \in \Ext} \in \ell^\infty (\Ext)\}
\ee
which is clearly a linear contraction, namely $ \|\Phi(f)\|_\infty
\le \|f\|$. On the other hand, for every $f\in\Aff$
$$\underset{\tau\in \text{co}(\Ext)}{\sup}|f(\tau)|\le \|\Phi(f)\|_\infty .$$
It follows from the well-known Krein-Milman theorem that
$$\|f\|=\underset{\tau\in\TA}\sup|f(\tau)|\le \|\Phi(f)\|_\infty.$$
Therefore,  $\Phi$ is actually an isometry, namely,
$\|f\|=\underset{\tau\in\TA}\sup|f(\tau)|=\|\Phi(f)\|_\infty.$

Furthermore,   when $\TA$ has only finitely many extremal points,
say $$\Ext= \{\tau_j\}_1^m,$$ i.e., when $\TA$ is a classical simplex,
then  $\Phi$ is onto $\mathbb R^m= \ell^\infty(\Ext)$.

\subsection{Strict comparison of projections and weak unperforation of $K_0(\A)$}\label{Ss:strict comp}

As in the literature  there is more than one definition  of
comparison of projections, we state explicitly   below the one that
we will use.

\bD {D: strict comp}   We say that a C*-algebra $\A$ has the strict comparison of projections if $\TA$ is non-empty and for any two
projections $p$ and $q$ in $\A$, the strict inequalities $\tau(p)<
\tau(q)$ for all $\tau\in \TA$ imply that $p \prec q$, namely
$[p] < [q]$, the strict ordering induced by the Murray-von Neumann
equivalence of projections.  \eD

The above definition was given in \cite[FCQ2,1.3.1]{Blackadar1988}
for simple C*-algebras, while we only work with simple C*-algebras
of real rank zero. For ease of reference, let us state without proof
the following simple observation.

 \bL{L:tau<1} Let $\A$ be a simple,
unital C*-algebra of real rank zero such that $\A$ has the strict
comparison of projections. If $p\in \A\otimes\K$ is a projection
such that $(\tau\otimes \tr)(p)< 1 $ for all $\tau \in \TA$, then there is a projection $p'\in \A$ such that
$p'\otimes e_{11}\sim p$ in $\A\otimes\K$. In particular,
$(\tau\otimes \tr)(p)=\tau(p')$ for all $\tau \in \TA$. \eL

As a trivial fact, if $\A\otimes \K$ has the strict comparison of
projections, then so does $\A$. Less obvious is the converse as
follows.

\bL{L:strict in stabilizer} Let $\A$ be a $\sigma$-unital simple
C*-algebra of real rank zero. Then the following are equivalent:
\item [(i)]  $\A$ has the strict comparison of projections.
\item  [(ii)] $M_n(\A)$ has the strict comparison of projections for every $n\in \mathbb N$.
\item [(iii)] $\A\otimes \K$ has the strict comparison of projections.
\eL \bp The implications (iii) $\Rightarrow$ (ii)$\Rightarrow$ (i)
are obvious. We need only to prove that (i) implies (iii).

To simplify notations, denote by $\tau\in \TA$ both a tracial weight on $\A \otimes \K$ as well as
its restriction to $\A$.  Let $p, q $ be projections in $\A \otimes \K$, such that $\tau (p)<\tau (q)$ for all $\tau\in \TA$.
Take an approximate identity $\{ e_\lambda\mid \lambda\in \Lambda \}$ of $\A$. Then \linebreak$\{ \sum_{i=1}^n e_\lambda \otimes e_{ii} \mid
(\lambda , n)\in \Lambda  \times \mathbb N \}$ is an approximate identity of $\A\otimes \K$. For $p, q
\in \A \otimes \K $, a  standard argument shows that $p, q$ are
equivalent to projections $p', q'$ in $ e_{\lambda _0}\A e_{\lambda
_0} \otimes \K$ for some $\lambda_0$. Therefore, from now on we can assume that $\A$ is unital.

Recall that the evaluation map  $\hat q-\hat p$ belongs to $ \Aff$ and choose
 $$0< \delta <\min \big( \inf_{\tau\in \TA } (\tau (q)-\tau (p)), 2\big).$$

Applying the diagonalization of projections in $\A $ or in $\A
\otimes \K$, respectively, proved in \cite[1.1,1.4]{Zhang2} and
\cite[1.2]{Zhang1},
  one can
 find  for every $k \in \mathbb N$ projections $p_k, q_k, r_k, r'_k$ in $\A \otimes
 \K$ such that
 \begin{align*}[p] &=2^k[p_k]+[r_k],\ \ 0< [r_k]<[p_k]\\
 [q]&=2^k[q_k]+[r'_k],\ \ 0< [r'_k]<[q_k].
\end{align*}
By iterating the diagonalization above to $r_k$ one can also assume
that $[r_k]< [r'_k]$.

 Choose $k$ such that $\frac 1{2^k} \sup_{\tau\in \TA } \tau (p) <\frac{\delta}{2} $ and
 $\frac 1{2^k} \sup_{\tau\in \TA } \tau (q) <\frac{\delta}{2} $. In particular, $\tau (p_k)\le  \frac{\delta}{2} <1$
 and $\tau (q_k)\le  \frac{\delta}{2} <1 $ and the same holds for $r_k$ and $r_k'$.
 By Lemma \ref{L:tau<1},   $p_k$, $q_k$. $r_k$, and $r_k'$ are equivalent to projections in $\A$. Thus, without loss of generality,
 we can assume  $p_k, q_k, r_k, r_k'\in \A$. Since $[r_k]<  [r'_k]$, it follows that $\tau(r_k)<\tau(r_k')$. Applying
 the strict comparison of projections of $\A$, one sees  $r_k\precsim r_k'$ in $\A$.

Moreover,  for every $\tau\in \TA$,
$$
\tau(q)-\tau(p)= 2^k\tau(q_k) +\tau(r'_k)- 2^k\tau(p_k) -\tau(r_k)> \delta$$
hence
 $$
2^k\big( \tau(q_k) - \tau(p_k))> \delta - \tau(r'_k)> \frac{\delta}{2}>0.
 $$
Thus $\tau(q_k)> \tau(p_k)$ for every $\tau\in \TA$. By the strict
comparison of projections of $\A$, it follows that $p_k\precsim q_k$ in
$\A$.  But then $p\precsim q$ in $\A\otimes \K$, which concludes the
proof. \ep

Finally, we notice the following direct consequence.

\bR{R: unperforated} Let $\A$ be a simple  $\sigma$-unital
C*-algebra with real rank zero. If $\TA$ has the strict comparison
of projections, then $K_0(\A)$ is weakly unperforated. \eR

\bp We sketch the reasoning only for the reader's convenience. For
any $x= [p]-[q] \in K_0(\A)$ with $nx>0$ for some $n\in \mathbb N$,
where projections $p, q\in \A\otimes \K$. Then $nx=[r]$ for some
projection $0\ne r\in \A\otimes \K$ by the strict comparison of
projections. We show that $x>0$. In fact, $n[p]=n[q]+[r]$ and for
all $\tau\in \TA$ (identified  with T$(\A\otimes \K)$). It follows
that $n\tau(p)=n\tau(q)+\tau(r)$ and hence
$\tau(p)=\tau(q)+\frac{1}{n}\tau(r)>\tau(q)$. By the strict
comparison of projections of $\A$, $q\precsim p$ and hence $x>0$. \ep

\bR{R:cancellation} For any unital C*-algebra the property of having
stable rank one implies cancellation (\cite[6.5.1]{BlackadarBook}).
Assume further that a C*-algebra has  real rank zero. Then the
cancellation property is equivalent to  stable rank one
(\cite[6.5.2]{BlackadarBook}).

\eR
\subsection{Quasitraces}\label{Ss:quasitraces}

Recall that a quasitrace $\tau$  on a C*-algebra $\A$ satisfies the same properties of a trace with exception of additivity,
 and that  $\tau(a+b)=\tau(a)+\tau(b)$ for $a, b\in \A_{sa}$ under the additional hypothesis that $a$ and $b$ commute.
A 2-quasitrace on $\A$ is a quasitrace that has a quasitrace extension to $\A \otimes \M_2(\mathbb C)$ and hence extends
 to $\M_n(\mathbb C)$ for every $n\in \mathbb N$ and thus to $\A\otimes \K$.  $\QA$ denotes the collection of 2-quasitraces
 on $\A$. $\QA$ too is a Choquet simplex and contains $\TA$ as closed face \cite [Proposition II 4.5]{BlackHandelm}.

Notice that for every projection $p\in \A\otimes \K$, the evaluation map  $$\QA\ni \tau \to \tau(p)$$ also belongs to
 $\QAff$ and we will still denote it by $\hat p$.

We will use the following density property \cite [Theorem 6.9.3]
{BlackadarBook} (see also \cite [Lemma III.3.4]{BlackHandelm}): if
$\A$ is  simple, unital, non-elementary (i.e., $\A\otimes\K
\not\cong \K$), stably finite, of real rank zero, stable rank one,
and has  weakly unperforated $K_0(\A)$, then $\{ \hat x   \mid x \in
K_0(\A )\}$ is uniformly dense in $ \QAff $, namely, for every
$f\in   \QAff$ and every $\epsilon
>0$ there is an element $x \in K_0(\A )$ such that \be
{e:density Q} |\tau(x)- f(\tau)|< \epsilon \ \  \forall \tau\in \QA.
 \ee

It is convenient to mention explicitly the following consequence of the above density property.

\bR{R:density for TA} If  $\A$ is  simple, unital, non-elementary,
stably finite, of real rank zero, stable rank one,  with  weakly
unperforated $K_0(\A)$, and with $\TA\ne \emptyset$, then the semigroup $
D(\A\otimes \K)$ of all equivalence classes of projections in $\A
\otimes \K $, which we identify with $\{ \hat p \mid [p] \in D(\A
\otimes \K )\}$,  is dense in $\Aff_+$; equivalently, for every
$f\in \Aff_+$ and every $\epsilon
>0$ there is a projection $p\in \A\otimes \K$ such that \be
{e:density T}
 |\tau(p)- f(\tau)|< \epsilon \ \  \forall \tau\in \TA. \ee
 \eR

\bp It is enough to notice that $\TA$ is a closed face of $\QA$
\cite[Proposition II 4.5]{BlackHandelm} and  $\QA$ is a Choquet
simplex. Then every $f\in \Aff$ has an extension to an $\tilde f\in
\QAff$ such that
 $\inf_{\tau\in \TA}f(\tau)\le \tilde f \le \sup_{\tau\in \TA}f(\tau)$ \cite[Theorem 11.22] {Goodearl}.
  In particular, if $f\in \Aff_+$,  then $\tilde f\in \QAff_+$. But then the same projection $p\in \A\otimes \K$ that satisfies
  (\ref {e:density Q}), satisfies a fortiori the condition in (\ref {e:density T}).
\ep A celebrated (unpublished) result of Haagerup \cite{UH} states
that if $\A$ is unital and exact, then $\QA=\TA$.  Then   Blanchard
and Kirchberg  \cite [Remark 2.29 (i)]{EBEK} observed that this
result can be extended non-unital exact C*-algebras. Brown and
Winter \cite {BrNWW} provided a short proof of Haagerup's result in
the finite nuclear dimension case.

We show now that $\QA=\TA$ holds also for the C*-algebras considered in the present paper. First, we recall  the
following fact which is an immediate consequence of  \cite [Theorem 11.22]{Goodearl}.
\bL{L: simplex}
Let $K$ be a Choquet simplex,  let $F\subset K$ be a closed face of $K$, and let $x\in \ext K)\setminus F$. Then
for all $\alpha, \beta \in \mathbb R_+$ there is a $g\in \text{Aff } (K)_+$ such that $g\mid_F=\alpha$ and  $g(x)=\beta$.
 \eL
\bT{T:TA=QTA} If $\A$ is a unital  simple,  C*-algebra of real rank
zero, stable rank one, and  has the strict comparison of
projections, then $\QA=\TA$. \eT

 \bp The case when $\A$ is elementary
algebra is trivial. We can assume henceforth that $\A$ is
non-elementary. Reasoning by contradiction, assume that $\TA$ is a
proper subset of $\QA $. Then there is an extreme point $\tau_o$
of $  \QA$ but not in $\TA$. Then $\{\tau_o\}$ and $\TA$ are closed
disjoint faces of $\QA$. Notice that constant functions on a face are
continuous and affine.

Thus by Lemma \ref {L: simplex}, there are positive continuous affine functions $f, g\in \QAff$ such that
$$\begin {cases} f\mid_{\TA}= \frac{1}{2} & f(\tau_o)=0\\
g\mid_{\TA}= 0 & g(\tau_o)=\frac{1}{2}.\end{cases}$$
By Remark \ref {R: unperforated}, $\A$ is also weakly unperforated, hence the conditions for the
 density of the projections in $\QAff_+$ (see (\ref {e:density Q})) are satisfied. Thus there are projections
  $p, q \in \A\otimes \K$ such that $$\begin {cases}\sup _{\tau \in \QA}|\tau(p)-f(\tau)|< \frac 1{4}\\
\sup _{\tau \in \QA}|\tau(q)-g(\tau)|< \frac 1{4}.
\end{cases}$$
In particular, $$\tau(q)< \frac 1{4}< \tau(p)\ \
\forall~\tau\in \TA.$$ By Lemma \ref {L:tau<1} there are projections
$p', q'\in \A$ such that $p'\sim p$, $q'\sim q$. Thus assume without
loss of generality that $p, q\in \A$. Then by the strict comparison
of projections, $q\precsim p$. This implies $\tau_o(q)\le \tau_o(p)$,
whereas $\tau_o (p)< \frac 1{4}< \tau_o (q)$, a contradiction. \ep

\section{Sums of commutators}\label{S:commutators}
We start with the following simple extension of a result by  Thomsen on uniform algebras.

\bL{L:Thomsen} Let $\C = \bigoplus_{i=1}^N p_iM_{n_i}(C(X_i))p_i$ where for each i, $n_i\in \mathbb N$, $X_i$
 is a compact Hausdorff space with covering dimension $d_i\le d$, and $ p_i \in \M_{n_i}(C(X_i))$ is a nonzero
   projection.  Let $a \in \A$  be a self-adjoint element, let $\eta>0$, and assume that
    $ |\tau(a)|\le \eta$ for all $ \tau \in T(\C)$. Then for every $\epsilon > 0$ there exist
$v_1, v_2, ..., v_d \in \C$ such that
$\| v_i \| \le \sqrt{2} \| a \|^{1/2}$  for $1 \le i \le d$ and
$$\| a - \sum_{i=1}^{d+1} [v_i, v_i^*] \| < \eta +\epsilon.$$
\eL
\bp
If  $\tau\in T(\C)$ and  $\tau(p_i)\ne 0$,  the restriction of $\frac{\tau}{\tau(p_i)}$ to
$p_iM_{n_i}(C(X_i))p_i$ is a tracial state on $p_iM_{n_i}(C(X_i))p_i$. Since $\sum _{i=1}^N \tau(p_i)=1$,
 $T(\C)$ is the collection of the convex combinations of the elements of $T(p_iM_{n_i}(C(X_i))p_i)$.
 Thus  it is enough to prove the statement for the case that $N=1$, i.e. when $\C = pM_n(C(X))p $.

Since  the range of the continuous function $\tr(p(x))$ consists of
integers, it follows that $X_0:=\{x\in X\mid \tr(p(x)) = 0\}$ is a
closed connected component of $X$ and $$p M_n(C(X))p= p
M_n(C(X\setminus X_0))p.$$ Thus we can assume without loss of
generality that $p(x)\ne 0$ for all $x\in X$ and hence $1\le
\tr(p(x)\le n$ for all $x \in X$. For every extremal $\tau\in
T(\C)$, there is an $x\in X$ for which $\tau(c)=
\frac{\tr(c(x))}{\tr (p(x))}$ for every $c\in \C$. Define $f(x) =
\frac{\tr(a(x))}{\tr(p(x))}$, then $f\in C(X)$ and $ | f(x) | \le
\eta$ for every $x\in X$. Define  $b:= fp$. Then $b\in \C$, $ \|
b\|\le \eta $, and $$\tr(b(x))= \tr(f(x)p(x))=f(x)\tr(p(x))=
\tr(a(x)).$$ Thus $\tr((a-b)(x))=0$ for every $x\in X$ and hence it
follows  from  Thomsen's result \cite [Lemma 1.4]
{ThomsenCommutator} that there exist $v_1, v_2, ..., v_d \in p
M_n(C(X))p$ such that $\| v_i \| \le \sqrt{2} \| a \|^{1/2}$ for $1
\le i \le d$ and $\| a -b- \sum_{i=1}^{d+1} [v_i, v_i^*] \| <
\epsilon$. As a consequence, $\| a - \sum_{i=1}^{d+1} [v_i, v_i^*]
\| <\eta + \epsilon.$ \ep

\bL{L:technical} Let $\C$ be the C*-inductive limit $\C = lim_{n
\rightarrow \infty} (\C_n, \phi_{n,n+1})$ of  a sequence of unital
C*-algebras $\{\C_n\}$ such that  $T(\C_n)\ne \emptyset$ for all $n$
and  the connecting maps $\phi_{n, n+1} : \C_n \rightarrow \C_{n+1}$
are unital and one-to-one.

\begin{enumerate}
\item [(i)]
Let $n\in \mathbb N$ and  $a\in \C_n$ be a self-adjoint element
such that
$\tau(a) >0$ for all $\tau \in T(\C)$.
Then there exists an $ M \geq n$ such that
$\tau(\phi_{n, m}(a)) > 0$ for all  integers $m \geq M$ and for all $\tau \in T(\C_m)$.
\item [(ii)] Let $\epsilon >0$, $n\in \mathbb N$, $a\in \C_n$ be a self-adjoint element
such that $|\tau(a) | < \epsilon$ for all $\tau \in T(\C)$.
Then there exists $M \geq n$ such that
$|\tau(\phi_{n, m}(a))| < \epsilon$
for all $\tau \in T(\C_m)$ and for all $m \ge  M$.
\end{enumerate}
\eL
\bp
To simplify notations, assume without loss of generality that the C*-algebras $\C_n$ form a
 nested sequence of subalgebras of $\C$ for which $\C= \overline{\cup_n \C_n}$ and thus the connecting
  map $\phi_{n, n+1}$ is just the identity embedding of $\C_n$ into $\C_{n+1}$.
\item [(i)]  Suppose, to the contrary, that for some increasing sequence of integers $n_k\ge n$ there is a
 sequence $\tau_k \in T(\C_{n_k})$  for which
\be{e:ineq}  \tau_k(a)  \le  0.
\ee
 Let $\tilde \tau_k$ be an arbitrary extension of $\tau_k$ to a  state   of $\C$.
 The state space of $\C$ being w*-compact, there is some w*-converging subnet  of the
 sequence $\tilde \tau_k$. To simplify notations, we can assume  that $\tilde \tau_k \to \tau$ where $\tau$ is a state of $\C$.

We claim that $\tau$ is a trace.
For every $a\in \C$, let $a_n\in \C_n$ be a sequence converging to $a$. Then
\begin{alignat*}{3}\tau(aa^*)&= \lim_n  \tau(a_na_n^*) & \text{(by the continuity of $\tau$)}\\
&=  \lim_n  \lim_k \tilde\tau_{k}(a_na_n^*) & \text{(by the definition of $\tau$)}\\
&=  \lim_n  \lim_k \tau_{k}(a_na_n^*) &\hspace{1cm} \text{(because $a_n\in \C_{n_k}$ for all $n_k\ge n$)}\\
&=  \lim_n  \lim_k \tau_{k}(a_n^*a_n) & \text{(because $\tau_{k}$ is a trace on $\C_{n_k}$)}\\
&= \tau(a^*a) & \text{(reversing the above argument.)}
\end{alignat*}
From the definition of $\tau$ and from (\ref{e:ineq}), we have that $\tau(a)\le 0$, a contradiction.
\item [(ii)]  For every $\tau\in T(\C)$,  $\tau(\epsilon I - a) >0$ and $\tau(\epsilon I + a) >0$. Thus by (i),
 there is an $M_1\ge n$ (resp. $M_2\ge n$) such that $\tau(\epsilon I - a) >0$  for every $\tau\in T(A_m)$ and every
  $m\ge M_1$, (resp. $\tau(\epsilon I + a) >0$ for every $\tau\in T(A_m)$ and every $m\ge M_2$). The conclusion
  then follows by taking $M= \max\{M_1, M_2\}$. \ep

From \cite [Theorems 2.9, 3.4]{CuntzPedersen} we know that every selfadjoint  element in the kernel of all tracial
states of a simple unital C*-algebra is the norm limit of sums of selfcommutators. We need however to obtain a bound
 on the number of the commutators and on the norms of the elements in the commutators. This was obtained for unital
 simple AH-algebra $\C$ with real rank zero and bounded
dimension growth and for some other algebras in \cite {MarcouxSmallNumberCommutators}. Based on the work in \cite {LinRR0AHEmbedding}, we can extend these results to a larger class of  C*-algebras. We start with an approximation property for elements with uniformly bounded traces, which include of course those in the kernel of all traces.

 \bL{lem:MarcouxPropositionReplacement}
 Let $\A$ be a unital separable simple C*-algebra
of real rank zero, stable rank one, and   the strict comparison of
projections.  Let $a \in \A$  be a self-adjoint element, let
$\eta>0$, and assume that  $ |\tau(a)|\le \eta$ for all $ \tau \in
T(\A)$. Then for every $\epsilon > 0$ there exist $v_1, v_2, v_3,
v_4 \in \A$ such that $\| v_i \| \le \sqrt{2} \| a \|^{1/2}$ for $1
\le i \le 4$ and
$$\| a - \sum_{i=1}^{4} [v_i, v_i^*] \| < \eta +\epsilon.$$
\eL

\bp
Assume without loss of generality that $\|a\|=1$. Since $\A$ has real rank zero, there is a selfadjoint element of
finite spectrum $a'$ with $\|a'- (1-\frac{\epsilon}{6})a\| <\frac{\epsilon}{6}$. Thus $\|a'\|\le 1$ and $\|a'- a\| <\frac{\epsilon}{3}$.

Notice that for every $\tau \in \TA$ we have
 $$|\tau(a')|\le \eta + |\tau(a)-\tau(a')|\le \eta + \|a-a'\|< \eta +\frac{\epsilon}{3}.$$
Write $a'= \sum_1^k\lambda_j p_j$ for some mutually orthogonal projections   $p_j\in \A$ and $\lambda_j\in \mathbb R$.

By  \cite [Theorem 4.5]{LinRR0AHEmbedding}, which, as pointed out by Emmanuel C. Germain in the review
MR1869626 (2002i:46053), holds also when $\A$ is not necessarily nuclear,
(see also \cite[Theorem 4.20] {ElliottGongRR0} and \cite {ElliottGongRR0 IV}), there exists a unital
 simple AH-algebra $\C$ with real rank zero and  dimension bounded by three, and a unital
*-embedding $\Psi : \C \rightarrow \A$ such that  $\Psi$ induces an isomorphism of the K-theory invariant:
$$K_*(\Psi) : (K_0(\C), K_0(\C)_+, K_1(\C), [1_{\C}])
\rightarrow (K_0(\A), K_0(\A)_+, K_1(\A), [1_{\A}]).$$
 Notice that the induced map on tracial simplexes, $T(\Psi) : \TA \rightarrow T(\C)$  is an affine homeomorphism
  of compact convex sets  (see \cite[Theorem 6.9.1] {BlackadarBook}) and all the tracial states of $T(\C)$
  extends
   to tracial states on $\A$. To simplify notations, assume that $\C$ is a subalgebra of $\A$ sharing the
   unit with $\A$ and that  $\Psi$ is the natural inclusion map.

Decompose $\C$ into a C*-inductive limit
$\C = lim_{n \rightarrow \infty} (\C_n, \phi_{n, n+1})$ where
each connecting map $\phi_{n, n+1} : \C_n \rightarrow \C_{n+1}$
is unital and injective, and where each $\C_n$ is a finite direct sum
of unital homogeneous C*-algebras with spectrum being a path-connected
finite CW-complex with dimension less than or equal to three.
Denote by $\phi_n$ the unital map embedding $\C_n$ in $\C= \overline{\cup_n \phi_n(\C_n)}$ and hence in $\A$.

Notice that by the assumption of stable rank one and hence  cancellation,
two Murray-von Neumann equivalent projections in $\A$ are
necessarily unitarily equivalent. Moreover, by the
 isomorphism of the ordered K-groups, every projection $p\in \A$ is equivalent to a projection $q\in \C$  and
 furthermore, every projection in $\C$ is unitarily equivalent to  a projection in $\phi_n(\C_n)$ for some
 $n\in \mathbb N$ (e.g., see \cite [Proposition 6.2.9, Corollary 5.1.7, and Appendix L. 2.2] {WeggeObook}.)
 Thus, by considering $p:=\sum_1^kp_j$ we can find  a unitary $u\in \A$ for which   $up_ju^*\in \phi_n(\C_n)$
 for all $1\le j\le k$. Let $q_j:=\phi_n^{-1}(up_ju^*)$ and  let $a'':= \sum_1^k\lambda_j q_j$. Since $a'= u\phi_n(a'')u^*$,
 $$|\tau(\phi_n(a''))| = |\tau(a')| <  \eta + \frac{\epsilon}{3} \quad  \text{for every } \tau\in \TA$$ and hence
$$|\tau(\phi_n(a''))|<   \eta + \frac{\epsilon}{3} \quad  \text{for every } \tau\in T(\C).$$
By Lemma \ref {L:technical}  there is some $m\ge n$ (actually, for every $m'\ge m$)  for which
$$|\tau(a'')|<   \eta + \frac{\epsilon}{3} \quad  \text{for every } \tau\in T(\C_m).$$ But then by Lemma \ref {L:Thomsen},
there are  four elements $v_i\in C_m $ with $$\|v_i\|\le \sqrt 2 \|a''\|^{\frac{1}{2}}= \sqrt 2 \|a'\|^{\frac{1}{2}}\le \sqrt 2
 \|a\|^{\frac{1}{2}}$$ and such that $\|a''- \sum_{1}^4 [v_i,v_i^*]\|\le  \eta + \frac{2\epsilon}{3}$.
Thus $$\|a'- u\sum_{1}^4 [\phi_m(v_i),\phi_m(v_i^*)]u^*\|\le\eta + \frac{2\epsilon}{3}$$  and
 hence $$\|a- u\sum_{1}^4 [\phi_m(v_i),\phi_m(v_i^*)]u^*\|\le\eta + \epsilon.$$

   \ep

The reduction argument in \cite {MarcouxSmallNumberCommutators}  shows that every element in the
kernel of the unique tracial state is the sum of two commutators and every selfadjoint element is the sum
 of four selfcommutators. The same reduction provides the following result in our setting.

  \bT{T:TraceZeroTwoCommutators}
   Let $\A$ be a unital simple separable C*-algebra
with real rank zero, stable rank one, and strict comparison of
projections. Let $a \in \A$ be an element such that $\tau(a) = 0$
for all $\tau \in \TA$. Then $a$ is the sum of two commutators;
i.e., there exist $y_1, y_2, y_3, y_4 \in \A$ such that
$$a = [y_1, y_2] + [y_3, y_4].$$
If, in addition, $a$ is self-adjoint then
$a$ can be expressed as the sum of four self-commutators; i.e.,
there exist $x_1, x_2, x_3, x_4 \in \A$ such that
$$a = \sum_{i=1}^{4} [x_i, x_i^*].$$

  \eT
\bp
The first step is to show that if $a\in \A$ is selfadjoint, then there are twelve elements $x_1, x_2, ..., x_{12} \in \A$ with
$\| x_i \| \leq  13 \| a \|^{1/2}$
such that
$$a = \sum_{i=1}^{12} [x_i, x_i^*]$$
 The proof is essentially the same as that of \cite{MarcouxSmallNumberCommutators}
3.9, which Marcoux presents as an adaptation of T. Fack  proof of
\cite[Theorem 3.1]{Fack} and its   modification by K. Thomsen \cite [Theorem 1.8]
{ThomsenCommutator}.  We just have to replace the  key step in
Marcoux's proof, namely \cite[Proposition 3.6]{MarcouxSmallNumberCommutators}
 obtained in the case of a unique tracial state, with
Lemma \ref{lem:MarcouxPropositionReplacement} that we have
obtained above. Since  the former approximates selfadjoint elements
by the sum of two selfcommutators while the latter needs four, the
same proof now decomposes $a$ into a sum of 12 selfcommutators.

The further reduction to two commutators or in the selfadjoint case
to four self-commutators follows then from
\cite[Theorem 3.10]{MarcouxSmallNumberCommutators}. \ep

While Marcoux's theorems on which we depend do not present explicitly  norm bounds for the elements composing the commutators, these bounds are implicit in his proofs and are further quoted explicitly in his Remark 5.3  in \cite {MarcouxSmallNumberCommutators}. We do not need a value for these bounds as their existence suffices for our needs.

\bR{R: norm bounds} There is a constant $M$, independent of the
algebra $\A$ and of the element $a\in \A$, such that the elements $y_i,
x_i\in \A$ in the above theorem satisfy $\|y_i\|\le M\|a\|^{1/2}$
and $\|x_i\|\le M\|a\|^{1/2}$ for $1\le i\le 4$. \eR

As a corollary, we get the following (see \cite{MarcouxIrishSurvey}
Section 2):

\bC{C:span comm is closed}
 Let $\A$ be a unital separable simple C*-algebra with
real rank zero, stable rank one, and strict comparison of
projections.
Then $[ \A, \A ]$ (the linear span of the commutators of $\A$) is
norm-closed. \eC

\section{Linear combination of projections}\label{S:lin comb proj}

In the previous section, we have seen that elements belonging to the
kernel of all tracial states are sums of commutators.  As shown by
Marcoux in \cite[Theorem 3.8]{MarcouxIndianaSpanProjections}, under  mild
conditions  every commutator is a linear combination of projections.
In fact, implicit
 in his proof and in the proof of his preceding lemmas is also an estimate on the number of projections
 needed and on the coefficients in that linear combination. Such an estimate is stated in \cite [Theorems 3.1, 3.3, 3.4]
 {MarcouxIrishSurvey}.  A discussion of that estimate is also given in our previous paper \cite {KNZPISpan}.
\bL {L:MarcouxCommutatorUniversalConstant}\cite[Lemma 2.4] {KNZPISpan} Let
$\A$ be a unital $C^*$-algebra for which there exist three mutually
  orthogonal projections  $p_1, p_2$ and $p_3$  such that $I  = p_1 + p_2 + p_3$ and
$p_i \precsim I  - p_i$ for $1 \leq i \leq 3$. Then for all $x, y \in
\A$ with $\| x \|, \| y \| \leq 1$,  there exist  $n\le 84$
projections $q_1, q_2, ..., q_n \in \A$ and  real numbers $\alpha_1,
\alpha_2, ..., \alpha_n$ with $|\alpha_j|\le 2\sqrt 2$ such that:
$$[x,y]:=xy-yx = \alpha_1 q_1 + \alpha_2 q_2 + ... + \alpha_n q_n.$$
\eL

Note that    such projections $p_1, p_2, p_3$ above exist in every
unital simple real rank zero $C^*$-algebra of dimension at least 3 by \cite [Theorem 1.1] {Zhang2}

  Thus combining Theorem \ref
{T:TraceZeroTwoCommutators}  and  \ref
{L:MarcouxCommutatorUniversalConstant} we obtain the following
result.

 \bT{T:Trace zero
LC proj}
 Let $\A$ be a unital simple separable C*-algebra
with real rank zero, stable rank one, and the strict comparison of
projections. Then every element  $a \in \A$  such that $\tau(a) = 0$
for all $\tau \in \TA$ is a linear combination $\sum
_{j=1}^{168}\alpha_jp_j$ of projections $p_j\in \A$ with
$|\alpha_j|\le 2 \sqrt{2} M^2 \|a\|$ where $M$ is the constant referred to
in Remark \ref {R: norm bounds}.

\eT

If $\A$ is unital and has a unique tracial state $\tau$, then every $a\in \A$ has
natural decomposition $$a= \tau(a)I+ (a-\tau(a)I) $$ into a  scalar multiple of a  projection and an element belonging to the kernel of the trace.
Under our  additional hypotheses on $\A$, the same holds also in the case when $\TA$ has only finitely many extremal points.
Furthermore, we can control the coefficients in the linear combination of projections.

\smallskip

\bL{L: m extremal} Let $\A$ be a simple  unital  C*-algebra of real
rank zero having strict comparison of projections and  assume that
$\TA$ has extremal points $\{\tau_1, \tau_2, \cdots, \tau_m\}$. Then for every
$\nu>0$ and every $a\in \A_{sa}$ there exist real numbers $\{\lambda_1, \lambda_2, \cdots, \lambda_m\}$
and projections
 $\{p_1, p_2, \cdots, p_m\}$ in $\A$ such that $\tau(a- \sum_{j=1}^m \lambda_j p_j)=0$ for every $\tau\in \TA$
 and $\sum_1^m |\lambda_j|\le (m+\nu) \|a\|$.
\eL
\begin{proof}
The case when $\A$ is a matrix algebra being trivial, assume without loss of generality that $\A$ is not elementary. Let $$
\Phi:  \Aff \ni f \to \Phi(f)=\{f(\tau)\}_{\tau \in \Ext} \in \ell^\infty (\Ext)= \mathbb R^m\}
$$ be the linear isometry defined in (\ref {e:isom}) and let $e_j$ be the standard basis of
$\mathbb R^m$. Then  $f_j:=\Phi^{-1}(e_j)$ is a basis of  $\Aff$.

Notice that $ f_j(\tau_i)=\delta_{i,j}\ge 0$ for every $1\le i,j\le m$, thus $f_j(\tau)\ge0$ for every $\tau\in \TA$ and every $j$. Thus by Remark
\ref {R:density for TA} for every $0< \delta< 1$    there exist projections $p_j\in \A\otimes \K$ such that
 $\|\hat p_j- (1-\delta)f_j\|\le \frac{\delta}{2} $ for every $j$. As a consequence,
  $\|\Phi(\hat p_j)- e_j\|_\infty < \frac{3\delta}{2} $.  Thus  $$\tau(p_j)\le (1-\delta)f_j(\tau)+
   \frac{\delta}{2} \le 1- \frac{\delta}{2}\quad \text{for every }\tau\in \TA\text{ and for every }j.$$
By Lemma \ref {L:tau<1}, there are projections in $\A$ with the same traces as $p_j$, thus we can
assume without loss of generality that $p_j\in \A$.

Let $B$ be the $m\times m$ matrix with columns $\Phi(\hat p_j)$.
 We can choose $\delta $ small enough so that $\|B-I\|\le \frac{\nu}{2m^{3/2}}$. Reducing if necessary $\nu$, we see that $B$
 is invertible and $\|B^{-1}-I\|\le \frac{\nu}{m^{3/2}}$. Thus $\{\Phi(\hat p_j)\}$ is a basis of $\mathbb R^m$ and
 hence $\{\hat p_j\}$ is a basis of $\Aff$. Thus, given a selfadjoint $a\in \A_{sa}$, there are (unique)
scalars $\lambda_j\in \mathbb R$ so that $\hat a= \sum_1^m
\lambda_j\hat p_j$, that is
$$ \tau(a) = \sum_1^m \lambda_j\tau( p_j) \quad \forall \tau\in \TA.
$$
Set $\lambda:=\begin{pmatrix}\lambda_1\\ \vdots\\\lambda_m\end{pmatrix}$. Hence
$$\Phi(\hat a) =\sum_1^m \lambda_j\Phi(\hat p_j) = B \lambda \quad \text{and thus} \quad  \lambda = B^{-1}\Phi(\hat a).$$

Thus
\begin{align*}
\sum_1^m |\lambda_j|&\le m\|\lambda\|_\infty \\ &=m\| B^{-1}\Phi(\hat a)\|_\infty\\
&\le m\| \Phi(\hat a)\|_\infty+ m\|(B^{-1}-I)\Phi(\hat
a)\|_\infty\\&\le
m\| \Phi(\hat a)\|_\infty+ m\|(B^{-1}-I)\Phi(\hat a)\|_2\\
&\le m\| \Phi(\hat a)\|_\infty+ m\|B^{-1}-I\|\,\|\Phi(\hat a)\|_2\\
&\le m\| \Phi(\hat a)\|_\infty+ m^{3/2}\|B^{-1}-I\|\,\|\Phi(\hat a)\|_\infty\\
&\le (m+\nu) \| \Phi(\hat a)\|_\infty\\&= (m+\nu) \| \hat a\|\\
&\le (m+\nu)\|a\|.
\end{align*}
\ep

By combining Theorem \ref {T:Trace zero LC proj} and Lemma \ref {L: m extremal} we thus obtain
that the C*-algebra $\A$ is the linear span of its projections and
furthermore  has a universal constant $V_o$ as in the following
theorem.

\bT{T:lin comb}
  Let $\A$ be a unital simple separable C*-algebra
with real rank zero, stable rank one, strict comparison of
projections and assume that $\TA$ has a finite number of extremal
points.   Then $\A$ is the linear span of its projections.
Furthermore, there exists a positive integer $N \ge 1$ and a
positive constant $V_0 > 0$ such that for every $a \in \A$  there
exist an integer $n \le N$, complex numbers $\alpha_1, \alpha_2,
...., \alpha_n$ with $| \alpha_i | \le V_0\|a\|$ for $1 \le i \le
n$, and projections $p_1, p_2, ..., p_n \in \A$ such that $a =
\alpha_1 p_1 + \alpha_2 p_2 + ... + \alpha_n p_n$. \eT

\section{C*-algebras that are not  the span of their projections}\label{S:notspan}
\subsection{Infinitely many extremal traces}\label{s:inf traces}
If we relax the condition that $\TA$ has finitely many extremal points,
$\A$ may fail to be the span of its projections.
\bP {P: no lin}  Let $\A$ be a simple $\sigma$-unital C*-algebra of
real rank zero  such that $\Ext$ is infinite, the collection $D(\A)$
of Murray-von Neumann equivalence classes of projections of $\A$ is
countable, and $D(\A\otimes \K)$ is dense in $\Aff_+$ (see Remark
\ref {R:density for TA}). Then $\A$ is not the linear span of its
projections. \eP
\begin{proof}
Let $\{[p_j]\}$ be an enumeration of $D(\A)$ and let $x_j:=\hat p_j$. Recall that $x_j$ does not depend on the representative projection $p_j$, that $x_j\in \Aff$,
 and that  $\|x_j\|\le 1$ for all $j$.

Recall that by \cite{Zhang3}
every projection $p\in \A\otimes \K$ is unitarily equivalent to
the direct sum of a finite number of projections $r_j\in \A$, or
more precisely, to a projection $\sum_{1}^n r_j\otimes e_{jj} $ with $r_j\in \A$ for every $j$.
Thus $\hat p =\sum_i^n \hat r_j$ and for each $j$, $\hat r_j = x_{j^\prime}$ for
some $j^\prime$. As a consequence, $\Aff$ is separable. Recall that
$\Aff$ is a real Banach space.

Now we follow the standard proof that the cardinality of any Hamel basis of an infinite dimensional separable
Banach space is not countable. Choose a unit length $y_1\in \{x_j\}$
such that $y_1\not\in \text{span}\{x_1\}$, and then choose $y_2\in
\{x_j\}\setminus  \text{span}\{x_1, x_2, y_1\}$. Recursively,
construct a sequence $y_k \in \{x_j\}$ such that $\|y_k\|=1$ and
$$y_k\not\in M_k:=\text{span}\{x_1, x_2, \cdots x_k,y_1, y_2\cdots, y_{k-1}\}.$$
Since $\Ext$ is infinite by hypothesis, $\Aff$ has infinite dimension too and hence also
  $\text{span}\{x_j\}$ has infinite dimension. Thus the construction cannot terminate after a finite
   number of steps and hence the sequence $\{y_k\}$ is infinite.

As $M_k$ is a closed subspace,  $\delta_k:= \text{dist}( y_k, M_k)>0$ for all $k$. Choose an increasing sequence
of integers $n_k$ such that $2^{n_{k+1}}> \frac{2^{n_{k}+1}}{\delta_k}$ and let $$y:=\sum_1^\infty \frac{y_j}{2^{n_j}}.$$
Then $y\in \Aff$.
Now $y_j= x_{\pi(j)}= \hat p_{\pi(j)}$ for some index $\pi(j)$. Set $a:=\sum_1^\infty \frac{p_{\pi(j)}}{2^{n_j}}$. Then $a\in \A_+$
and  for every $\tau\in \TA$ we have
$$\tau(a)=\sum _1^\infty \tau(\frac{p_{\pi(j)} } {2^{n_j} })= \sum _1^\infty \frac{\hat p_{\pi(j)} (\tau) }{2^{n_j} } = \sum _1^\infty
 \frac{y_j(\tau)} {2^{n_j}}
= y(\tau).$$ Reasoning by contradiction, assume  that $a$ is a linear combination  of projections $r_i\in \A$, say $a= \sum _{i=1}^m\lambda_ir_i$.
 Since $a= \sum _{i=1}^m\frac{\lambda_i+\bar \lambda_i}{2} r_i$, assume that $\lambda_j\in \mathbb R$.  Then for every $\tau\in \TA$
 we would also have $$y(\tau)= \tau(a) = \sum _{i=1}^m\lambda_i \tau( r_i).$$
For each $i$, $r_i\in [p_{i'}]$ for some $i'$ and hence $\tau (r_i)=\hat p_{i'}(\tau) = x_{i'}(\tau) $. Thus
 $$y= \sum _{i=1}^m\lambda_i x_{i'}  \in \text{span} \{x_j\}.$$ But then $y\in M_k$ for some $k$. Since
 $\sum_{j=1}^{k-1}\frac{y_j}{2^{n_j}}\in M_k$ by the definition of $M_k$, it follows that
\begin{align*} \frac{\delta_k}{2^{n_k}}&
=\text{dist}( \frac{y_k}{2^{n_k}}, M_k)\\
&= \text{dist}(y -  \sum_{j=1}^{k-1}\frac{y_j}{2^{n_j}}- \sum_{j=k+1}^\infty \frac{y_j}{2^{n_j}}, M_k)\\
&= \text{dist}(-\sum_{j=k+1}^\infty \frac{y_j}{2^{n_j}}, M_k)\\
& \le \| \sum_{j=k+1}^\infty\frac{y_j}{2^{n_j}}\|\\
& \le \sum_{j=k+1}^\infty\frac{1}{2^{n_j}}\\
& \le \frac{2}{2^{n_{k+1}}}
\end{align*}
a contradiction.
\end{proof}

Recall from Remark \ref {R:density for TA} that if $\A$ is  simple,
unital, non-elementary, of real rank zero, stable rank one, with
strict comparison of projections (which implies both that  $K_0(\A)$
is  weakly unperforated and that $\TA\ne \emptyset$), then
$D(\A\otimes \K)$ is dense in $\Aff_+$.   Examples of  C*-algebras
satisfying the conditions in Proposition \ref {P: no lin} can be
found in the category of simple unital AF-algebras. Indeed, by a
result of Blackadar \cite {Blackadar1980}, every Choquet simplex can
be realized as the tracial state space of some simple unital
AF-algebra. Thus it is enough to start with  a Choquet simplex with
infinitely many extremal points, e.g., the Brauer simplex with
extreme boundary $[0,1]$.    Explicit examples can be found among
crossed products coming from Cantor minimal systems \cite
{GiordanoPutnamSkau}.

\bR{R:answ Marcoux}
Real rank zero C*-algebras satisfying the conditions of Proposition \ref{P: no lin} provide a negative answer to Marcoux's questions
 \cite[Question 1 and 2] {MarcouxIrishSurvey} on whether the span of the projections
 in  a simple unital C*-algebra must be closed and on whether if the span is dense, it
  must coincide with the algebra.
\eR

\subsection{Not LP but with ``many projections"}\label{s: nor RR0}
C*-algebras  can fail to be the span of their projections even if they contain a ``many projections". 

A C*-algebra has the LP property if the span of its projections is dense. Of course, real rank zero algebras have the LP property and projectionless algebras do not. But it may be interesting to note that there are algebras with ``many projections", e.g, having the same ordered $K_0$ group as a real rank zero algebra and satisfying the SP property (every hereditary subalgebra contains a nonzero projection) and yet fail to satisfy the LP property.

\bR {R:NoLP}
Let $\A$ be a C*-algebra with two distinct tracial states $\tau, \tau'$  such that $\tau(p) = \tau'(p)$ for every projection $p \in \A$.
Then $\A$ does not have the LP property. \eR

\bp
Since   $\tau$ and $ \tau'$ agree on the projections of $\A$ and hence on their linear combinations,
they must agree also on the norm closure of span of the projections of $\A$.  Hence the latter cannot coincide with $\A$.
\ep

An example of a ``nice" algebra having two distinct tracial states
that agree on all the projections of the algebra is an
\emph{AI}-algebra $\A$, (an inductive limit of finite direct sums of
the form $\mathbb{M}_{n_1}(C[0,1]) \oplus \mathbb{M}_{n_2}(C[0,1])
\oplus ... \oplus \mathbb{M}_{n_k}(C[0,1])$) constructed by
\cite{ThomsenAI} (a special case of \cite{Villadsen}) which is is
simple, unital, with $K_0(\A) = \mathbb{Q}$ (rational numbers),
$K_0(\A)_+ = \mathbb{Q}_+$ (positive rationals), $[I_{\A}] = 1 \in
\mathbb{Q}$ and has tracial simplex $\TA \cong [0,1]$ having two extreme
points $\tau, \tau'$. It follows  that $\tau(p) = \tau'(p)$ for
every projection $p \in \A$. It can also be shown that $\A$ has the SP property. 

\subsection{Non-unital algebras}\label{s:Non-unital}
If we relax the condition that $\A$ is unital, we also see that $\A$
may fail to be the span of it projections. The simplest example is provided by the algebra $\K$ where infinite rank operator clearly cannot be a linear combination of projections in $\K$, which are finite.

More generally,  no simple stable $\sigma$-unital C*-algebra of real rank zero with non-empty tracial simplex $\TA$ can be the span
of its projections. To see this, first recall that F. Combes  showed
in an early work \cite [Proposition 4.1 and Proposition 4.4]
{Combes} that every semifinite (also called densely defined)  lower
semicontinuous weight $\tau$ on a C*-algebra $\A$ has an extension
to a normal weight $\bar \tau$ on the enveloping von Neumann algebra
$\A^{**}$ and  if the weight is tracial, then the extension is
unique. More recently, Ortega, Rordam, and Thiel proved in
\cite[Proposition 5.2]{OrtegaRordamThiel} that if the weight $\tau$
is tracial then the extension $\bar \tau$ is also tracial. Notice
that the faithfulness of $\tau$ does not guarantee the faithfulness
 of $\bar\tau$. However, if $\tau$ is faithful,  $Q\in \A^{**}$ is an open projection, and $\bar\tau(Q)=0$, then $Q=0$.

Notice also that for every element $a\in \A$, the range projection $R_a$ of $a$ as an element in the enveloping von Neumann
 algebra  $\A^{**}$ is an open projection. Indeed $$R_a= \chi_{aa^*}(0, \|a\|^2]= \chi_{aa^*}(0, \|a\|^2+1)$$ where $\chi_{aa^*}$
 denotes the spectral measure of $aa^*$ in $\A^{**}$.  $R_a$ may  fail to belong to $\Mul(A)$ (actually, most likely it is not).

\bL{L: no comb} Let  $\A$  be a $\sigma$-unital simple C*-algebra of
real rank zero. No element $a\in \A $ with
 $\bar \tau (R_a) =\infty $ for at least one  $\tau\in \TA$  can be a linear combination of
 projections. Such an element $a$ always exists in $\A \otimes \K$ when $\TA \not= \emptyset$.
\eL
\begin{proof}
Assume  that $a=\sum_{i=1}^m \lambda _i p_i$ for some scalars
$\lambda _i$ and projections $p_i\in \A $. Then $R_a\le
\bigvee_{i=1}^mp_i$ where the supremum is of course taken in
$\A^{**}$. Then
$$\bar\tau(R_a)\le \bar\tau\Big(\bigvee_{i=1}^mp_i\Big)\le \sum_{i=1}^m\bar\tau(p_i)= \sum_{i=1}^m\tau(p_i)< \infty$$
where the second inequality is a well-known von Neumann algebra property derived from the Kaplanski parallelogram law.

To see the last statement, first one sees that a sub C*-algebras of
$\A \otimes \K$ is *-isomorphic to $\K$ by Brown's Stabilization
Theorem (\cite{BrownStable}). Then  $a = \sum_1^\infty
\frac{1}{n}e_{ii}$ satisfies $\bar \tau (R_a) =\infty$  for all
$\tau\in \TA$.
\end{proof}

\subsection{Non-simple algebras}\label{s:non-simple}
Again, the simplest example of ``nice" non-simple algebras that fail to be the linear span of their projections is given by $\K$, or more precisely by the unitization $\A:=\mathbb C I+\K$ of $\K$. This has been observed by Marcoux in \cite {MarcouxIrishSurvey}. Indeed, the collection of linear combinations of projections in $\A$ is $\mathbb C I+\mathscr F$ where $\mathscr F$ denotes the finite rank class and hence is a proper subset of  $\A$.

\section{Positive combinations of projections}

For each positive element $a$ in a C*-algebra $\A$ of
real rank zero, the hereditary C*-subalgebra $\her(a) := (a\A a)^-$
is also of real rank zero by \cite{BrownPedersen}, and hence, has a
sequential approximate identity of projections $\{p_i\}$.  In the enveloping von Neumann algebra $\A^{**}$, the sequence $\{p_i\}$
converges strongly to the range projection $R_a$ of  $a$.  Furthermore, $$\her(a)= \her (R_a): = (R_a\A^{**} R_a)^-\cap \A. $$

Given a nonzero open projection $Q\in \A^{**}$ with the property that $\bar\tau(Q)< \infty$ for
every $\tau\in \TA$, the evaluation map $\TA\ni \tau\to \hat Q(\tau):= \bar\tau(Q)$ is clearly affine.
 In the case that $\Ext$ is finite, every affine map is also continuous and hence $\hat Q\in \Aff$. Since $Q$ is a
 nonzero open projection and every trace $\tau\in \TA$ is faithful, it follows that $\inf_{\tau\in \TA}\hat Q>0$.
  To summarize, when $\Ext$ is finite and $Q\in \A^{**}$ is an open projection, then

\be{e:inf} \bar\tau(Q)< \infty ~\forall\, \tau\in \TA~\Rightarrow~\begin {cases}\sup_{\tau\in \TA}\hat Q(\tau)<\infty\\\inf_{\tau\in
\TA}\hat Q(\tau)>0\end{cases}
 \ee

\bT{T:main}
 Let $\A$ be a $\sigma$-unital simple separable C*-algebra
with real rank zero, stable rank one, strict comparison of
projections and assume that $\TA$ has a finite number of extremal
points. Then an element $a\in \A_+ $ is a positive combination of
projections if and only
 if $ \bar \tau (R_a) <\infty$ for all  $\tau\in \TA$.
\eT

 We will present the proof through the chain of the following lemmas. Our first result extends to the present setting the main tool that we used in   \cite {KNZPISpan} and \cite {KNZFiniteSumsVNA}.

\bL{L:invert} Let  $\A$ be a simple C*-algebra $\A$ of real rank
zero and stable rank one having strict comparison of projections and
such that $\TA$ has
 a finite set of extremal tracial states. Let  $p, q$ be projections in $\A$ with $qp=0$, $q\precsim p$ and let  $b=qb=b q$ be
 a positive element of $\A$. Then for every  scalar $\alpha
>\|b\|$, the positive element $a:=\alpha p\oplus b$ is a positive combination of projections.
\eL

\begin{proof}
Let $r:=p+q$, then the corner $r\A r$ of $\A$
satisfies the same hypotheses as $\A$ and it is unital. Thus by   Theorem \ref {T:lin comb} there is a universal constant $V_0$ such that
for every  $a \in r\A \,r$,
there exist scalars $\alpha_1, \alpha_2, ..., \alpha_n$ and
 projections $p_1, p_2, ..., p_n \in r\A r$ such that
$$a = \alpha_1 p_1 + \alpha_2 p_2 + ... + \alpha_n p_n$$
and
$$| \alpha_1 | + | \alpha_2 | + ... + | \alpha_n | \leq V_0 \| a \|.$$
This is precisely condition (1) of \cite[Proposition 2.7]{KNZPISpan}. Condition (2) of the same proposition,
namely that the positive combinations of projections of $r\A r$ are norm dense in $(r\A r)_+$, is an immediate
consequence of the hypothesis that $RR(\A)=0$ and hence $RR(r\A r)=0$. Thus the conclusion of
\cite[Proposition 2.7]{KNZPISpan} applies, namely every positive invertible operator in $r \A r $ is a positive
 combination of projections.  This permits to apply  \cite[Lemma 2.9]{KNZPISpan} which yields the requested
 positive combination of projections.
\end{proof}

\bL{L: open} Let $\A$ be a stable $\sigma$-unital C*-algebra of real rank zero, $\emptyset \ne \TA $ has the
 strict comparison of projections, and let   $Q\in \A^{**}$ be an open projection.
\item [(i)]  $\sup _{\tau \in \TA} \bar \tau (Q)<\infty$   if and
 only if there exists a projection $r\in \A$ such that $\bar \tau (Q)< \tau(r)$ for all $\tau\in \TA$.
 \item [(ii)] Condition (i) is satisfied if and only if there is a projection $r\in \A$, an open projection $R'\le r$,
  and a partial  isometry $w\in \A^{**}$ with $w w^*=Q$ and
 $w^*w =R'$ such that the map $\phi (x) = w^*xw$  is a trace-preserving *-isomorphism between $\her(Q) $ and $\her(R')$.
 \item [(iii)] If $\Ext$ is finite and condition (ii) is satisfied, then the projection $r$ can be chosen so that
 $\bar\tau(Q)< \tau(r) < 2\bar\tau(Q)$ for all $\tau\in \TA$.
\eL

 \begin{proof} \item (i) The sufficiency is obvious since    $\bar \tau (Q)< \tau(r)=\hat r(\tau)$ implies that
 $\bar \tau (Q)$ is bounded since $\hat r$ is a continuous function on the compact set $\TA$.

Now we prove the necessity.  The stability of $\A$ guarantees that for every $\tau\in \TA$ the trace $\bar \tau$ is
infinite, that is $\bar\tau(I)=\infty$.
 Let $\{p_i\}$ be an approximate identity of $\A$ consisting of projections. Then $p_j \uparrow I$ and
 hence $\tau(p_j)\uparrow \infty$. Thus $\{\hat p_i\}$ is a monotone increasing
 sequence of continuous functions on $\TA $ with $\lim\hat p_i(\tau )=\infty $ pointwise. Since $\TA$ is compact,
 it follows from the well-known Dini's Theorem in elementary topology that there exists an integer $n_0$
 such that $$\inf _{\tau \in \TA}   \tau (p_{n_0}) >\sup _{\tau \in \TA} \bar \tau (Q). $$ Thus it is enough to
  set $r=p_{n_0}$.
\item [(ii)] If $Q\sim R'\le r$, then $\bar\tau(Q)= \bar\tau(R')\le \tau(r)$ for all $\tau\in \TA$. In order to
 obtain strict inequality, it is enough to replace $r$ with a projection $r'\in \A$, $r'\gneqq r$. This proves the sufficiency.

For the necessity, assume that $\sup _{\tau \in \TA} \bar \tau (Q)<\infty$ and let $r\in \A$ be the projection with
 $\bar \tau (Q)< \tau(r)$ for all $\tau\in \TA$ provided by (i).
  Since the hereditary algebra $\her(Q)=(Q\A^{**}Q)\cap \A$ has real
 rank zero \cite {BrownPedersen} and is $\sigma$-unital, one can find an
 increasing approximate identity of projections $q_i\in \A$ for $\her(Q)$ and setting $r_i=q_{i+1}-q_i$ obtain $Q =\bigoplus_{i=1}^\infty
  r_i$ and hence
  $$
\sum _{i=1}^\infty \tau(r_i)= \bar \tau (Q)\le \sup _{\tau \in \TA} \bar \tau (Q)
  .$$

Since $\tau (r_1)<\tau (r)$ for all   $\tau \in \TA $, by the strict
comparison of projections one can find a partial isometry $v_1\in
\A$ such that $v_1v_1^* =
  r_1$ and $v_1^*v_1 =r'_1<r $. Similarly, because of $\tau (r_2)<\tau (r-r'_1)$ for all $\tau \in \TA $, one
  can find another partial isometry $v_2\in \A$ such that $v_2v_2^* =
  r_2$ and $v_2^*v_2 =r'_2<r-r_1 $. Repeating   the construction recursively, one obtains a sequence of partial isometries
   $\{v_i\}\subset \A  $ with  mutually
  orthogonal initial projections $\{r_i\}$ in $\her(Q)$ and mutually orthogonal range projections $\{r'_i\}$  in $r\A r$.
  Define $R':= \sum_{i=1}^\infty r'_i$ and $w=\sum_{i=1}^\infty v_i$. Then $R'\in \A^{**}$ is an open projection, $R'\le r$,
   $w\in \A ^{**}$  is a partial isometry, and  $ww^*=Q$ and $w^*w= R'$. Thus $Q\sim R'$.
Define $\phi: \her (Q) \rightarrow \A^{**}$ by $\phi (x) = w^* xw$.
Using the fact that $q_j=\sum_{i=1}^j r_i$ converges in the strict topology to $Q$, it is now routine to show that
 $\phi(x)\in \A$ for every $x\in \her (Q)$ and hence that $\phi $ is a trace-preserving *-isomorphism from $\her (Q)$ onto $\her(R')$.
\item [(iii)]
By the hypothesis that $\Ext$ is finite, it follows that the evaluation map $\hat Q(\tau)=\bar\tau(Q)$ is continuous
on $\TA$. Identify $\A$ with $\A\otimes \K$ and choose $\{q_j\}$ to be an approximate identity of $\her (Q\oplus Q)$
consisting of projections of $\A\otimes\K  $. Then $\{\tau (q_j)\}$
increases to $\bar \tau (Q\oplus Q) = 2\bar \tau (Q)$, that is the sequence of continuous functions $\hat q_j$ increase
 pointwise to the continuous function $2\hat Q$. By Dini's theorem, the convergence is uniform. Thus choose  $r_0=q_{j_0}$
 for an appropriate $j_0$, as wanted.
\end{proof}

The next lemma is the technical crux of the proof.

\bL{L:unital and cond}   Let  $\A$ be a unital simple C*-algebra
$\A$ of real rank zero and stable rank one with
strict comparison of projections, and with finitely many extremal tracial states.
 Let $a\in \A_+ $ be such that
$  \bar \tau (R_a) > \frac{1}{2}$ for all  $\tau\in \TA$. Then $a$ is a positive combination of projections.
\eL
 \begin{proof}
Assume  without loss of generality that $\|a\|=1$.
Let $\{\tau_i\}_1^m$ be the collection of the extremal tracial
states of $A$.
Notice that
\begin{equation}\label{e:rR}
I= \chi_{[0,1]}(a)= \chi_{\{0\}}(a)+R_a.
\end{equation}
By the w*-continuity of each $\bar\tau_i$,
$$\lim_{\lambda\to0+}\bar\tau(\chi_{(0, \lambda)}(a))=0\quad\text{and}\quad \lim_{\lambda\to0+}\bar\tau(\chi_{( \lambda, 1]}(a))
=\bar\tau(R_a).$$ Since by (\ref{e:rR})
$$\bar\tau_i(\chi_{\{0\}}(a))=1- \bar\tau_i(R_a)< \bar\tau_i(R_a)
\quad \forall ~1\le i\le m,$$ we can find $0< \alpha < \beta<1$ such
that
$$
\bar\tau_i(\chi_{[0, \alpha)}(a))< \bar\tau_i(\chi_{(\beta,
1]}(a))\quad \forall ~1\le i\le m
$$
and hence
\begin{equation}\label{e:alphabeta}
\bar\tau(\chi_{[0, \alpha)}(a))< \bar\tau(\chi_{(\beta, 1]}(a))\quad
\forall ~\tau\in \TA.
\end{equation}
Choose arbitrary numbers $0 < \gamma_1< \gamma_2< \gamma_3< \alpha< \gamma_4 < \beta$.
Let $f : [0, 1] \rightarrow [0, 1]$ be the  continuous
function defined by
$$
f(t) =
\begin{cases}
t & t\in [0,1]\setminus [\gamma_1, \gamma_3]\\
\gamma_1 &  t \in[\gamma_1, \gamma_ 2]\\
\text{linear }  & t \in[\gamma_2, \gamma _3].
\end{cases}
$$
Notice that  $f(a)\ge 0$ and  $R_{f(a)}=R_a$.

Recall that for every $0\le  \gamma< \delta< \infty $, the projection $\chi_{[\gamma, \delta]}(a)$ is closed. Now by using the hypothesis that   $\A$ is unital, it follows that $\chi_{[\gamma, \delta]}(a)$ is compact.
 Notice also that   $\chi_{(\gamma, \delta)}(a)$ is open and so is $\chi_{(\gamma,1]}(a)=\chi_{(\gamma,2)}(a)$.
Thus by Brown's interpolation property \cite{BrownIMP}, there exist
projections $s, p, q \in  \A $ such that
\begin{alignat*}{5}
&\chi_{[0, \gamma_1]}(a)&~\le~ &s&~\le~& \chi_{[0, \gamma_2)}(a)\\
&\chi_{[0, \gamma_3]}(a)&~\le~ &q&~\le~ &\chi_{[0, \alpha)}(a)\\
&\chi_{[\beta, 1]}(a)&~\le~ &p&~\le~& \chi_{(\gamma_4, 1]}(a).
\end{alignat*}
It is immediate to see that
\begin{equation}\label{e:spq} pq=0, \quad s\le q,\quad\text{and}\quad p\le I-s.
\end{equation}
Define the following elements of $\A$:
\begin{align*}
b~&:=a - f(a) + s f(a) s\\
a_1&:=b+ \alpha p\\
a_2 &: = (I-s)f(a)(I-s)-\alpha p.
\end{align*}
Since the spectral projections of $a$ commute with $a$ and hence with $f(a)$, and since $s- \chi_{[0, \gamma_1]}(a)\le \chi_{(\gamma_1, \gamma_2)}(a)$, it follows that $s- \chi_{[0, \gamma_1]}(a)$ commutes with both $f(a)\big(\chi_{[0, \gamma_1]}(a)+\chi_{[\gamma_2,
1]}(a)\big)$ and with $ \gamma_1\chi_{(\gamma_1, \gamma_2)}(a)$.  Hence it commutes with $$f(a) = f(a)\big(\chi_{[0, \gamma_1]}(a)+\chi_{[\gamma_2,
1]}(a)\big)+ \gamma_1\chi_{(\gamma_1, \gamma_2)}(a).$$
But then  it follows that also
$s$ and hence $I-s$ commute with $f(a)$. As a consequence we have that
$$b= a-f(a)+f(a)s\quad\text{and}\quad a_2= f(a)-f(a)s-\alpha p
$$ and hence
$$a= a_1+a_2.$$
Now  $b\ge 0$ because $a \ge f(a)\ge 0$. Since $$a-f(a)=(a-f(a))\chi_{[\gamma_1,\gamma_3]}(a)\le\chi_{[\gamma_1,\gamma_3]}(a)\le q$$ and
$$f(a)s\le s\le q,$$
hence $b=qbq$.
Moreover
\begin{align*}\|a-f(a)\|&= \gamma_2-\gamma_1\\
\|f(a)s\|&\le \|f(a)\chi_{[0, \gamma_2]}(a)\|\le \gamma_1,
\end{align*}
and hence $\|b\|\le  \gamma_2< \alpha.$
Now   by (\ref {e:alphabeta}) we have
$$
\tau(q)= \bar\tau(q) \le \bar\tau(\chi_{[0, \alpha)}(a)) <
\bar\tau(\chi_{[\beta, 1]}(a))\le  \tau(p) \quad \forall ~\tau\in
\TA.
$$
By the strict comparison of projections  we obtain that $q\precsim
p.$

As $a_1= b+\alpha p$ with $b=qbq\ge0$, $\|b\|< \alpha$,  $qp=0$ and
$q\precsim p $, we obtain by Lemma \ref {L:invert} that
$a_1$ is a positive combination of projections in $\A$.

We prove now that the same holds for $a_2$. Notice first that by (\ref {e:spq})
\begin{equation}\label {e:range a2} R_{a_2}\le  I-s\in \A. \end{equation}

Since
$$\chi_{[0, \gamma_2]}(a)-s\le \chi_{[0, \gamma_2]}(a)-\chi_{[0, \gamma_1]}(a)= \chi_{(\gamma_1, \gamma_2]}(a)$$ and since
$f(a)\chi_{(\gamma_1, \gamma_2]}(a)=\gamma_1\chi_{(\gamma_1, \gamma_2]}(a)$, it follows that
$$
f(a) \big(\chi_{[0, \gamma_2]}(a)-s\big)= \gamma_1\big(\chi_{[0,
\gamma_2]}(a)-s\big).
$$ Thus
\begin{align*}
a_2&= f(a)(I-s)-\alpha p\\
&\ge f(a)(I-s)-\alpha \chi_{(\gamma_4, 1]}(a)\\
&=f(a)\big(\chi_{[0, \gamma_2]}(a)-s\big)
+ f(a)\chi_{(\gamma_2,\gamma_4]}(a)+ (f(a)-\alpha) \chi_{(\gamma_4, 1] }(a)\\
&= \gamma_1\big(\chi_{[0, \gamma_2]}(a)-s\big)
+ f(a)\chi_{(\gamma_2,\gamma_4]}(a)+ (f(a)-\alpha) \chi_{(\gamma_4, 1] }(a)\\
&\ge \gamma_1\big(\chi_{[0, \gamma_2]}(a)- s\big)+ \gamma_2\chi_{(\gamma_2,\gamma_4]}(a)+ (\gamma_4-\alpha)\chi_{(\gamma_4, 1]}(a)\\
&\ge\min\{ \gamma_1, \gamma_4-\alpha\}\Big(\chi_{[0, \gamma_2]}(a)- s+\chi_{(\gamma_2,\gamma_4]}(a)+\chi_{(\gamma_4, 1]}(a)\Big)\\
&= \min\{ \gamma_1, \gamma_4-\alpha\}(I-s).
\end{align*}
Thus $R_{a_2}\ge I-s$ and by (\ref {e:spq}) it follows that  $R_{a_2}= I-s$  and $$
a_2\ge  \min\{ \gamma_1, \gamma_4-\alpha\}R_{a_2},$$
i.e.,  $a_2$ is locally invertible. But then by  \cite [Lemma 2.9] {KNZPISpan} (see also proof of Lemma \ref {L:invert}) $a_2$
 is a positive combination of projections in $\A $. This concludes the proof.
\end{proof}

 By using these lemmas we can now provide the proof of Theorem \ref {T:main}.
\bp [Proof of the theorem]
By Lemma \ref {L: no comb} we see that the condition that $\bar\tau(R_a)<\infty$ for all $\tau\in \TA$ is necessary.

To prove the sufficiency, assume first that $\A$ is stable. Then by Lemma \ref {L: open} (iii) there is a projection
 $r\in \A$ such that $\bar\tau(R_a)<\tau(r)< 2\bar\tau(R_a)$ for all $\tau\in \TA$ and by part (ii) of the same lemma,
  and in the notations of the lemma, there is a trace preserving  *-isomorphism $\phi$ from $\her(a)= \her (R_a)$ onto
  $\her (R') =\her (\phi(a))$ with $R'\le r$. Decomposing $\phi(a)$ into a positive combination of projections
  (necessarily in $\her (\phi(a)$) is equivalent to decomposing $a$ into a positive combination of projections (necessarily
in $\her(a)$). Thus we can thus assume without loss of generality that $R_a\le r$. By passing to the corner $r\A r$ of $\A$,
which satisfies the same properties as $\A$, we can further assume that $\A$ is unital, i.e., identify $r$ with $I$. By renormalizing
the trace $\tau$, we thus have that  $\bar\tau(R_a)>\frac{1}{2}$. But
then the conclusion that $a$ is a positive combination of projections follows from Lemma \ref {L:unital and cond}.

Finally, we remove the condition that $\A$ is stable. Let  $a\in \A_+$ satisfy the condition $\bar\tau(R_a)< \infty$ for all
$\tau\in \TA$. Since $R_{a}$ is an open projection,  $R_{a}=\sum_1^\infty r_j$ for some projections $r_j\in \A$ and the series
converges in the strict topology. Then
\ba
(\overline{\tau \otimes \tr})\big(R_{ a \otimes e_{11} }\big)&=
(\overline{\tau \otimes \tr})\big(R_a \otimes e_{11} \big)
=(\overline{\tau \otimes \tr})\big(\sum_1^\infty r_j \otimes e_{11}\big)\\
&=\sum_1^\infty (\overline{\tau \otimes \tr})(r_j \otimes e_{11})= \sum_1^\infty (\tau \otimes \tr)(r_j \otimes e_{11})\\&= \sum_1^\infty \tau (r_j)=\bar\tau\big(\sum_1^\infty r_j\big)
=\bar\tau(R_a).
 \end{align*}

Since $\A \otimes \K$ satisfies the same conditions as $\A$ and is stable, it follows from the first part of the proof that
$a \otimes e_{11}$ is a positive combination of projections belonging to $\her (a\otimes e_{11})= \her (a)\otimes e_{11}$.
But then $a$ too is a positive combination of projections in $\her(a)$.
\ep

If $\A$ is unital, for every $a\in \A_+$ it follows that   $ \bar \tau (R_a) \le \bar \tau (I) =1$ for all  $\tau\in \TA$.  Thus $a $ is a positive linear combination of projections by  Theorem \ref {T:main}.
\bC{C:unital} Let $\A$ be a unital simple separable
C*-algebra with real rank zero, stable rank one, strict comparison
of projections and assume that $\TA$ has a finite number of extremal
points. Then every element $a\in \A_+ $ is a positive linear combination of projections in $\A$.
\eC
\bC{C:non pos} Let $\A$ be a $\sigma$-unital simple separable
C*-algebra with real rank zero, stable rank one, strict comparison
of projections and assume that $\TA$ has a finite number of extremal
points. Then an element $a\in \A $ is a linear combination of
projections if and only
 if $ \bar \tau (R_a) <\infty$ for all  $\tau\in \TA$.
\eC
\bp
The necessity is given by Lemma \ref {L: no comb}. Assume that $ \bar \tau (R_a) <\infty$ for all
$\tau\in \TA$ and let $Q:=R_{aa^*+a^*a}$. Since $Q= R_a\bigvee R_{a^*}$ and since $R_a\sim R_{a^*}$ in
$\A^{**}$, it follows that $\bar \tau(Q)\le 2\bar\tau(R_a)<  \infty$ for every $\tau\in \TA$. Now $a$
decomposes naturally into a linear combination of positive elements of $\A$, all with range projections
dominated by $Q$ and hence all satisfying the condition of Theorem \ref {T:main}. Thus they are all positive
combinations of projections and hence $a$ is a linear combination of projections.
\ep

We would like to point out why for unital algebras we first proved that every element is a linear combination of projections and from that we deduced that every positive element is a positive combination of projections, while in the non-unital case we first proved that positive decompositions hold for positive elements. In other words, why employing directly Theorem \ref {T:lin comb} for the proof of Corollary \ref{C:non pos} would  not suffice. Indeed while even in the non-unital case there is still a trace-preserving *-isomorphism $\phi$  from $\her( aa^*+a^*a)$ onto a
(hereditary) subalgebra of a corner $r\A r$ for some projection $r\in \A$, we would obtain from Theorem
\ref {T:lin comb} only that $\phi(a)$ is a linear combinations of projections in $r\A r$ but then we could not guarantee that those projections belong to $\phi\big(\her( aa^*+a^*a)\big)$  and hence that they can be carried back to $\her( aa^*+a^*a)$ and provide a
 decomposition for $a$ itself. 

\end{document}